\documentclass[12pt]{article}

\usepackage{amssymb}
\usepackage{amsmath}
\usepackage{lscape}
\usepackage{graphicx}
\usepackage{url}
\usepackage{fullpage}
\expandafter\ifx\csname color\endcsname\relax
    \newcommand{\rline}[1]{\overline{#1}}
    \newcommand{\gline}[1]{\underline{#1}}
    \newcommand{\rgline}[1]{\rline{\gline{#1}}}
    \newcommand{\bline}[1]{#1|}
    \newcommand{\pline}[1]{|#1}
    \newcommand{\bpline}[1]{\bline{\pline{#1}}}
    \newcommand{\bbline}[1]{#1||}
    \newcommand{\ppline}[1]{||#1}
\else
    \definecolor{olive}{named}{OliveGreen}	
    \definecolor{mahogany}{named}{Mahogany}	
    \newcommand{\rline}[1]{{\color{red}\overline{{\color{black}#1}}}}
    \newcommand{\gline}[1]{{\color{olive}\underline{{\color{black}#1}}}}
    \newcommand{\rgline}[1]{\rline{\gline{#1}}}
    \newcommand{\bline}[1]{#1{\color{blue}|}}
    \newcommand{\pline}[1]{{\color{mahogany}|}#1}
    \newcommand{\bpline}[1]{\bline{\pline{#1}}}
    \newcommand{\bbline}[1]{#1{\color{blue}||}}
    \newcommand{\ppline}[1]{{\color{mahogany}||}#1}
\fi

\input xy
\xyoption{all}

\newcommand{\matX}{\mathbf{X}}
\newcommand{\matM}{\mathbf{M}}
\newcommand{\XX}{\mathcal{X}}
\newcommand{\MM}{\mathcal{M}}
\newcommand{\jalg}{\mathbf{H}_3(\mathbb{O})}
\newcommand{\Btz}{B^2_{tz} - B^3_{tz}}
\newcommand{\dotBtz}{\dot B^2_{tz} - \dot B^3_{tz}}
\newcommand{\myendofproof}{\hspace{4.65in}$\square$}

\newcommand{\RR}{\mathbb{R}}
\newcommand{\CC}{\mathbb{C}}
\newcommand{\HH}{\mathbb{H}}
\newcommand{\OO}{\mathbb{O}}
\newcommand{\ZZ}{\mathbb{Z}}

\newcommand{\ba}{\mathfrak{b}}
\newcommand{\ga}{\mathfrak{g}}
\newcommand{\ha}{\mathfrak{h}}
\newcommand{\ma}{\mathfrak{m}}
\newcommand{\pa}{\mathfrak{p}}
\newcommand{\ra}{\mathfrak{r}}
\newcommand{\ta}{\mathfrak{t}}

\renewcommand{\aa}{\mathfrak{a}}
\newcommand{\bb}{\mathfrak{b}}
\newcommand{\cc}{\mathfrak{c}}
\newcommand{\dd}{\mathfrak{d}}
\newcommand{\ee}{\mathfrak{e}}
\newcommand{\ff}{\mathfrak{f}}
\newcommand{\gtwo}{\mathfrak{g}_2}

\newcommand{\sla}{\mathfrak{sl}}
\newcommand{\so}{\mathfrak{so}}
\newcommand{\su}{\mathfrak{su}}
\newcommand{\uone}{\mathfrak{u}(1)}
\newcommand{\uonem}{\mathfrak{u}(-1)}

\newcommand{\suH}{\su(2)_\HH}
\newcommand{\suII}{\su(2)_2}

\newcommand{\phit}{\phi_s}
\newcommand{\phiType}{\phi_t}
\newcommand{\phiHp}{\phi_\HH}

\newcommand{\mystrut}{\hbox{\vrule height 2.8ex width 0pt}}

\begin{document}


\title{\boldmath\textbf{%
  Discovering Real Lie Subalgebras of $\ee_6$ using Cartan Decompositions}}

\author{%
Aaron Wangberg\\
Department of Mathematics \& Statistics\\
Winona State University\\
Winona, MN 55987\\
\texttt{awangberg@winona.edu}
\and
Tevian Dray\\
Department of Mathematics\\
Oregon State University\\
Corvallis, OR  97331\\
\texttt{tevian@math.oregonstate.edu}
}

\date{\today}

\maketitle

\begin{abstract}
The process of complexification is used to classify a Lie algebra and identify
its Cartan subalgebra.  However, this method does not distinguish between real
forms of a complex Lie algebra, which can differ in signature.  In this paper,
we show how Cartan decompositions of a complexified Lie algebra can be
combined with information from the Killing form to identify real forms of a
given Lie algebra.  We apply this technique to $\sla(3,\OO)$, a real form of
$\ee_6$ with signature (52,26), thereby identifying chains of real subalgebras
and their corresponding Cartan subalgebras within $\ee_6$.  Motivated by an
explicit construction of $\sla(3,\OO)$, we then construct an abelian group of
order~8 which acts on the real forms of $\ee_6$, leading to the identification
of 8 particular copies of the 5 real forms of $\ee_6$, which can be
distinguished by their relationship to the original copy of $\sla(3,\OO)$.
\end{abstract}

\newpage
\section{Introduction}

The group $E_6$ has a long history of applications in physics, beginning with
the original discovery of the Albert algebra and its relationship to
exceptional quantum mechanics~\hbox{\cite{Jordan,JNW,Albert}}.  As a candidate
gauge group for a Grand Unified Theory, $E_6$ is the natural next step in the
progression $SU(5)$, $SO(10)$, both of which are known to lead to interesting,
albeit ultimately unphysical, models of fundamental particles (see
e.g.~\cite{Georgi}), and is closely related to the Standard Model gauge group
$SU(3)\times SU(2)\times U(1)$.

A description of the group $E_{6(-26)}$ as $SL(3,\OO)$ was given
in~\cite{Denver}, generalizing the interpretation of $SL(2,\OO)$ as (the
double cover of) $SO(9,1)$ discussed in~\cite{Lorentz}.  An interpretation
combining spinor and vector representations of the Lorentz group in 10
spacetime dimensions was described in~\cite{York}, and in~\cite{E6struct} we
obtained nested chains of subgroups of $SL(3,\OO)$ that respect this
Lorentzian structure.

The resulting action of $E_6$ appears to permit an interpretation in terms of
electroweak interactions on leptons, suggesting that this approach may lead to
models of fundamental particles with physically relevant properties.  In
particular, the asymmetric nature of the octonionic multiplication table
appears to lead naturally~\cite{Spin} to precisely three generations of
particles, with single-helicity neutrinos, observed properties of nature which
as yet have no theoretical foundation.

In the present work, we use various Cartan decompositions of the Lie algebra
$\ee_6$ to further extend this construction in two distinct ways.  First, we
identify additional real subalgebras of~$\sla(3,\OO)$, thus completing the
explicit construction of nested chains of subgroups of $SL(3,\OO)$ begun
in~\cite{E6struct}, building on the Lorentzian structure of $SL(2,\OO)$.  In
particular, we locate the ``missing'' $C_4$ subgroups of $SL(3,\OO)$ referred
to in~\cite{E6struct}, in the form $SU(3,1,\HH)$.  We then reinterpret the
Cartan decompositions used in our construction as (vector space) isomorphisms
of the complexified Lie algebra $\ee_6^\CC$, yielding an abelian group of such
isomorphisms that acts on the 5 real forms of $\ee_6$, allowing us to identify
8 particular copies of these real forms, based on their relationship to the
original copy of $\sla(3,\OO)$.

We briefly review the basic structure of $SL(3,\OO)$ in Section~\ref{sl30},
and the basic properties of Cartan decompositions in Section~\ref{graded}.  We
construct three important Cartan decompositions of $\ee_6$ in
Section~\ref{Auto3}, and use them to construct nested chains of subalgebras of
$\sla(3,\OO)$ in Section~\ref{Auto3Apps}, which are used in
Section~\ref{Chains} to construct a tower of subalgebras in~$\sla(3,\OO)$ based
upon our preferred basis of the Cartan subalgebra, completing the program
begun in~\cite{E6struct}.  Finally, we explore the relationships between the
real forms of $\ee_6$ in Section~\ref{discussion}.

\section{The Basics}
\label{basics}

\subsection{\boldmath$SL(3,\OO)$}
\label{sl30}

\begin{table}[tbp]
\begin{center}
\begin{tabular}{|l|ccc|}
\hline
           & \hspace{1cm} $(2)$ \hspace{1cm} & \hspace{1cm} $(3)$
	     \hspace{1cm} & \hspace{1cm} $(21)$ \hspace{1cm} \\
\cline{2-4}
\mystrut
Boosts     & $\dot B^1_{tz}$ & $\dot B^1_{tx}$ & $\dot B^1_{tq}$ \\[0.5ex]
           & $\dot B^2_{tz}$ & $\dot B^2_{tx}$ & $\dot B^2_{tq}$ \\[0.5ex]
           &                 & $\dot B^3_{tx}$ & $\dot B^3_{tq}$ \\[0.5ex]
\hline
           &    $(7)$        &     $(3)$       &    $(21)$ \\
\cline{2-4}
\mystrut
Simple     & $\dot R^1_{xq}$ & $\dot R^1_{xz}$ & $\dot R^1_{zq}$ \\[0.5ex]
Rotations  &                 & $\dot R^2_{xz}$ & $\dot R^2_{zq}$ \\[0.5ex]
           &                 & $\dot R^3_{xz}$ & $\dot R^3_{zq}$ \\[0.5ex]
\hline
Transverse &   $(7)$         &     $(7)$       &   $(7)$ \\
\cline{2-4}
\mystrut
Rotations  &  $\dot A_q$     &    $\dot G_q$   & $\dot S^1_q$ \\[0.5ex]
\hline
\end{tabular}
\caption{Our preferred basis for $\sla(3,\OO)$}
\label{Basis}
\end{center}
\end{table}

We summarize here the description of the Lie group $SL(3,\OO)=E_{6(-26)}$, as
given in~\cite{Denver,York,E6struct}.%
\footnote{Much of this material also appears in~\cite{aaron_thesis}.}
The \textit{Albert algebra} $\jalg$ consists of the $3\times3$ octonionic
Hermitian matrices.  A \textit{complex} matrix $\MM$ is one whose components
lie in some complex subalgebra of $\OO$; such matrices act on $\XX\in\jalg$ as
\begin{equation}
\XX \longmapsto \MM\XX\MM^\dagger
\label{MXM}
\end{equation}
and $SL(3,\OO)$ can be defined as the (composition of) such transformations
that preserve the determinant of $\XX$.  Each \textit{Jordan matrix} $\XX$ can
be decomposed as
\begin{equation}
\XX = \left( \begin{array}{c|c}
	\matX & \theta \\ \hline \theta^\dagger & \cdot
	\end{array}\right)
\end{equation}
and $\matM\in SL(2,\OO)\subset SL(3,\OO)$ acts on \textit{vectors} $\matX$ and
\textit{spinors} $\theta$ via the embedding
\begin{equation}
\matM \longmapsto
\MM = \left( \begin{array}{c|c}
	\matM & 0 \\ \hline 0 & 1
	\end{array}\right)
\label{typeI}
\end{equation}
An explicit identification of elements of $SL(2,\OO)$ with the (double cover
of the) Lorentz group $SO(9,1)$ was given in~\cite{Lorentz}, naturally
generalizing the identification of $SL(2,\CC)$ with (the double cover of)
$SO(3,1)$.  We refer to elements of $SL(3,\OO)$ of the form~(\ref{typeI}) as
being of \textit{type}~1, with analogous definitions of types 2 and 3.  There
are thus 3 natural copies of $SL(2,\OO)$ sitting inside $SL(3,\OO)$, which we
label according to their type and character as a Lorentz transformation.
Thus, the 45 elements of type 1 $SL(2,\OO)$ consist of 9 \textit{boosts}
$B^1_{tq}$ and 36 independent \textit{rotations} $R^1_{qr}$, where the spatial
labels $q$, $r$ run over $x$, $z$, and the imaginary octonionic units
$\{i,j,k,k\ell,j\ell,i\ell,\ell\}$.  Some of these 45 \textit{generators} of
type 1 $SL(2,\OO)$ require the use of more than one transformation of the
form~(\ref{MXM}), a phenomenon we call \textit{nesting}.

The generators of the 3 natural copies of $SL(2,\OO)$ do indeed span
$SL(3,\OO)$, but are not independent, so we introduce a preferred set of
generators.  First of all, the 3 copies of $G_2\subset SO(7)\subset SO(9,1)$
are in fact the same, so we dispense with the superscript labeling type, and
call the 14 generators $\{A_q,G_q\}$, with $q$ now ranging over the imaginary
octonionic units.  We denote the remaining 7 generators of (type~1) $SO(7)$ by
$S^1_q$.  Second, by triality there is only one copy of $SO(8)$, so we add the
7 generators $R^1_{xq}$, for a total of 28 generators so far.  The remaining
24 rotations in $SL(3,\OO)$ are the 3 types of rotations with $z$, to which
must be added the $27-1=26$ independent boosts --- the diagonal $zt$ boosts
are not all independent.

Our 78 preferred generators of $SL(3,\OO)$ become a basis of the Lie algebra
$\sla(3,\OO)$ under differentiation, denoted by a dot; this preferred basis is
summarized in Table~\ref{Basis}.  For further details,
see~\cite{Denver,York,E6struct}.

\subsection{Graded Lie Algebras}
\label{graded}

A \textit{$\ZZ_2$-grading} of a Lie algebra $\ga$ is a decomposition
\begin{equation}
\ga = \pa \oplus \ma
\label{decomp}
\end{equation}
such that
\begin{eqnarray}
\left[ \pa,\pa \right] &\subset& \pa \nonumber\\
\left[ \ma,\ma \right] &\subset& \pa \label{grade}\\
\left[ \pa,\ma \right] &\subset& \ma \nonumber
\end{eqnarray}
so that $\pa$ is itself a Lie algebra, but $\ma$ is merely a vector space.  We
will assume further that $\ma$ is not abelian (in which case it would be an
ideal of $\ga$), and that the Killing form is nondegenerate on each of $\pa$
and~$\ma$.%
\footnote{For an example of a $\ZZ_2$-grading that does not satisfy these
conditions, see~\cite{E6struct}.}
Each $\ZZ_2$-grading defines a map $\theta$ on $\ga$, given by
\begin{equation}
\theta(P+M) = P-M \qquad (P\in\pa, M\in\ma)
\end{equation}
which is clearly an involution, that is, a Lie algebra automorphism whose
square is the identity map.  Conversely, an involution $\theta$ defines a
$\ZZ_2$-grading in terms of the eigenspaces of $\theta$ with eigenvalues
$\pm1$, which can be shown to satisfy~(\ref{grade})~\cite{gilmore}.

A \textit{Cartan decomposition} of a real Lie algebra $\ga$ is a
$\ZZ_2$-grading of $\ga$, with the further property that the Killing form $B$
is negative definite on $\pa$ and positive definite on $\ma$.  In this case,
($\pa$,$\ma$) is called a \textit{Cartan pair}, the signature of $\ga$ is
($|\pa|$,$|\ma|$), and the associated involution $\theta$ is called a
\textit{Cartan involution}.  Informally, $\pa$ consists of \textit{rotations},
and $\ma$ of \textit{boosts}.

We extend this terminology to the complexification $\ga^\CC$ of $\ga$.
Whereas each real form of $\ga^\CC$ admits a unique (up to isomorphism) Cartan
decomposition, $\ga^\CC$ itself will admit several, one for each real form.

A slight modification of an involution $\theta$ on $\ga$ can be used to map
one real form of~$\ga^\CC$ to another.  Use $\theta$ to split $\ga$ into
eigenspaces $\pa$ and $\ma$ as above, and then introduce the
\textit{associated Cartan map} $\phi^*$ on $\ga^\CC$ via
\begin{equation}
\phi^*(P+M) =  P + \xi M \qquad (P\in\pa^\CC, M\in\ma^\CC)
\label{associated}
\end{equation}
with $\xi^2=-1$, that is, where $\xi$ is a square root of~$-1$ which commutes
with all imaginary units used in the representation of~$\ga$.  The structure
constants of a real form $\ga$ of $\ga^\CC$ are real by definition, and since
\begin{equation}
\left[ \pa, \pa \right] \subset \pa
  \hspace{1cm}
  \left[ \pa, \xi \ma \right] \subset \xi \ma
  \hspace{1cm}
  \left[ \xi \ma, \xi \ma, \right] \subset \xi^2 \pa = (-1)\pa
\end{equation}
then $\pa \oplus \xi \ma$ also has real structure constants, and is therefore
also a real form of the complex Lie algebra~$\ga^\CC$.  Note that $\phi^*$ is
a vector space isomorphism, but not a Lie algebra isomorphism --- as must be
the case if it takes one real form to another.  Further information regarding
the interplay between involutive automorphisms, the Killing form, and real
forms of a complex Lie algebra may be found in~\cite{gilmore}.

We claim that \textit{any} $\ZZ_2$-grading of a real Lie algebra $\ga$ is in
fact the Cartan map associated with (the restriction of) \textit{some} Cartan
decomposition of $\ga^\CC$.

\textbf{Lemma:}
Every $\ZZ_2$-grading of a real Lie algebra $\ga$ is the restriction of a
Cartan decomposition of $\ga^\CC$.  Equivalently, the extension of the
$\ZZ_2$-grading to $\ga^\CC$ is the image of a Cartan decomposition of some
(other) real form of $\ga^\CC$ under an associated Cartan map.

\textbf{Proof:}
Let $\ga=\pa\oplus\ma$ be a $\ZZ_2$-grading of $\ga$.  Extend this grading to
the $\ZZ_2$-grading $\ga^\CC=\pa^\CC\oplus\ma^\CC$ of $\ga^\CC$, then restrict
each component to the compact real form $\ga_c$ of $\ga^\CC$.  We thus obtain
the vector space decomposition
\begin{equation}
\ga_c = \pa_c \oplus \ma_c
\label{Cdecomp}
\end{equation}
Since we started with a $\ZZ_2$-grading, and since $\ga_c$ is a Lie algebra by
assumption, we must have
\begin{equation}
[ \pa_c, \pa_c ] \subset \pa^\CC \cap \ga_c = \pa_c
\end{equation}
with similar expressions holding for the other commutators in~(\ref{grade}).
Thus,~(\ref{Cdecomp}) is a $\ZZ_2$-grading of $\ga_c$.  It is now
straightforward to invert the associated Cartan map, constructing the vector
space
\begin{equation}
\ga' = \pa_c \oplus \xi\ma_c
\label{Rdecomp}
\end{equation}
with $\xi^2=-1$.  Even though associated Cartan maps are not Lie algebra
isomorphisms, they do preserve the $\ZZ_2$-grading, so $\ga'$ is a Lie
algebra, and the Killing form is negative definite on $\pa_c$, and positive
definite on $\ma_c$, by construction (and the assumed nondegeneracy of the
Killing form).  Thus,~(\ref{Rdecomp}) is a Cartan decomposition of the real
form $\ga'$ of $\ga^\CC$.

\myendofproof

The $\ZZ_2$-gradings of $\ga$ therefore correspond to the possible real forms
of $\ga^\CC$.  We will use associated Cartan maps~$\phi^*:\ga^\CC\to\ga^\CC$
not only to identify different real forms of~$\ga^\CC$, but also to identify
real subalgebras of our particular real form~$\ga$.  When applied to the
compact real form $\ga_c=\pa_c\oplus\ma_c$, the map~$\phi^*$ changes the
signature from~$(|\pa|+|\ma|, 0)$ to~$(|\pa|,|\ma|)$.  Similar counting
arguments give the signature of~$\phi^*(\ga)$ when~$\ga$ is non-compact, since
$\phi^*$ changes some compact generators into non-compact generators, and
vice-versa.  We then use the rank, dimension, and signature of~$\phi^*(\ga)$
to identify the particular real form of the algebra.  Using tables of real
forms of~$\ga^\CC$ showing their maximal compact subalgebras~\cite{gilmore},
we can also identify~$\pa=\ga\cap\phi^*(\ga)$ as a subalgebra of our original
real form~$\ga$.  Since the maximal compact subalgebra of~$\phi^*(\ga)$ is
known, we can further identify its pre-image as a non-compact subalgebra
of~$\ga$.

\section{\boldmath Cartan Decompositions of $\sla(3,\OO)$}
\label{Auto3}

\subsection{\boldmath Some Gradings of $\sla(3,\OO)$}

We first make some comments about our preferred basis for~$\sla(3,\OO)$, which
is listed in Table~\ref{Basis} and further discussed in~\cite{E6struct}.
Let~$\ba$ and~$\ra$ be the vector subspaces consisting of boosts and
rotations, respectively.  Our preferred choice of basis favors type~$1$
transformations, in the sense that we choose to represent $\so(8)$, as well as
its subalgebra $\gtwo$, in terms of transformations of type~$1$.  Let~$\ta_1$
be the subspace spanned by $\dotBtz$ and all type~$1$ transformations, and
$\ta_2$ and $\ta_3$ the subspaces spanned by the type~$2$ and type~$3$
transformations in our preferred basis which are not in~$\ta_1$.%
\footnote{Since $\dot{B}^1_{tz}+\dot{B}^2_{tz}+\dot{B}^3_{tz}=0$,
and since $\ta_1$ contains $\dot{B}^1_{tz}$ and~$\dotBtz$, we regard both
$\dot{B}^2_{tz}$ and $\dot{B}^3_{tz}$ as being elements of $\ta_1$.}
Let $\ta_{23}=\ta_2\oplus\ta_3$.  Finally, let $\ha$ be the subspace
corresponding to transformations that preserve the preferred quaternionic
subalgebra $\HH=\langle 1,k,k\ell,\ell \rangle$, and let $\ha^\perp$ be its
orthogonal complement.  The vector space $\ha$ is spanned by transformations
with no labels in $\lbrace i, j, j\ell, i\ell \rbrace$ and is in fact a
subalgebra, while~$\ha^\perp$ is the subspace spanned by transformations
having one label from~$\lbrace i, j, j\ell, i\ell \rbrace$; $\ha$ contains 
\textit{quaternionic} basis elements, such as~$\dot A_{k\ell}, \dot
R^1_{z\ell}$, and~$\dot B^3_{tx}$, while~$\ha^\perp$ contains
\textit{orthogonal-quaternionic} basis elements such as~$\dot A_i$, $\dot
R^1_{zi\ell}$, and~$\dot B^3_{tj\ell}$.

Direct computation shows that:
\begin{equation}
\begin{array}{ccc}
\left[\ra, \ra\right] \subset \ra &
	\left[\ta_1, \ta_1 \right] \subset \ta_1 &
	\left[\ha, \ha\right] \subset \ha \\
\left[\ra, \ba\right] \subset \ba &
	\hspace{1cm}
	  \left[\ta_1, \ta_{23} \right] \subset \ta_{23}
	  \hspace{1cm} &
	\left[\ha, \ha^\perp\right] \subset \ha^\perp \\
\left[\ba, \ba\right] \subset \ra &
	\left[\ta_{23}, \ta_{23} \right] \subset \ta_1 &
	\left[\ha^\perp, \ha^\perp\right] \subset \ha \\
\end{array}
\end{equation}
so that each of the decompositions $\ra \oplus \ba$, $\ta_1 \oplus \ta_{23}$,
and $\ha \oplus \ha^\perp$ is a $\ZZ_2$-grading of $\sla(3,\OO)$.  We therefore
introduce the involutions $\phit$, $\phiType$, and $\phiHp$ on $\sla(3,\OO)$,
given by
\begin{equation}
\begin{array}{ccc}
\phit(R + B) = R - B &
	\hspace{1cm} &
	R \in \ra, B \in \ba \\
\phiType(T_1 + T_{23}) = T_1 - T_{23} & &
	T_1 \in \ta_1, T_{23} \in \ta_{23} \\
\phiHp(H + H^\prime) = H - H^\prime & &
	H \in \ha, H^\prime \in \ha^\perp\\
\end{array}
\end{equation}
and let $\phit^*$, $\phiType^*$, and $\phiHp^*$ be the associated
Cartan maps on $\ee_6^\CC$.

%



The associated Cartan map $\phit^*$ transforms $\sla(3,\OO)$, which has
signature $(52,26)$, into the compact real form
$\phit^*\left(\sla(3,\OO)\right)$, which has signature $(78,0)$.  The
subalgebra
\begin{equation}
\phit^*\left(\sla(3,\OO)\right) \cap \sla(3,\OO) = \ra
\end{equation}
has dimension~$52$ and is therefore easily seen to be the compact real
form $\su(3,\OO)$ of~$\ff_4$.  For this case, the compact part
of~$\phit^*\left(\sla(3,\OO)\right)$ is the entire algebra, so that its
pre-image is already a known subalgebra of $\sla(3,\OO)$, namely $\sla(3,\OO)$
itself.

\subsection{\boldmath Some Subalgebras of $\sla(3,\OO)$}

We obtain two interesting subalgebras of $\sla(3,\OO)$ when we apply the
associated Cartan map~$\phiType^*$.  First, the signature of $\ga' =
\phiType^*\left(\sla(3,\OO)\right)$ is again $(52,26)$, since $\ta_{23}$
contains the same number of boosts and rotations~(16 of each).  The real form
$\ga'$ is therefore isomorphic to, but distinct from, $\sla(3,\OO)$, with
maximal compact subalgebra $\ga'_c=\ff_4$.  Hence, the pre-image
of~$\ga'_c$ is a real form of~$\ff_4$ in~$\sla(3,\OO)$.  It has signature
$(36,16)$, and the~$16$ non-compact generators identify this real form of
$\ff_4$ as $\su(2,1,\OO)=\ff_{4(-20)}$.  The second subalgebra comes from
looking at
\begin{equation}
\phiType^*\left(\sla(3,\OO)\right) \cap \sla(3,\OO) = \ta_1
\end{equation}
Because~$|\ta_1| = 46$, it must be a real form of $\dd_5\oplus\dd_1$.  But
$\ta_1$ contains~$\so(9,1)=\sla(2,\OO)$, which has signature~$(36,9)$, as well
as the boost $\dotBtz$.  Writing $\uonem$ for the non-compact real
representation of $\dd_1$ generated by $\dotBtz$, we identify~$\ta_1$ as the
subalgebra $\sla(2,\OO)\oplus\uonem$, with signature~$(36,9)\oplus(0,1)$.

We also obtain two subalgebras by applying $\phiHp^*$ to
$\sla(3,\OO)$.  In this case, the signature of
$\ga' = \phiHp^*\left(\sla(3,\OO)\right)$ is $(36,42)$, since
$\ha^\perp$ contains 28 rotations and 12 boosts, resulting in a net of 16
changes in signature.  The maximal compact subalgebra $\ga'_c$
of~$\ga'$ is~$\su(4,\HH)$, the compact real form of~$\cc_4$.  The pre-image
of $\ga'_c$ has signature $(24,12)$, and, due to the~$12$ non-compact
generators, can be identified as $\su(3,1,\HH)$.  The invariant subalgebra
\begin{equation}
\phiHp^*\left(\sla(3,\OO)\right) \cap \sla(3,\OO) = \ha
\end{equation}
has dimension $|\ha|=38$ and signature $(24,14)$, and is therefore a real form
of $\aa_5 \oplus \aa_1$ with signature $(21,14) \oplus (3,0)$.  Of the $24$
rotations unaffected by $\phiHp^*$, there are $21$ which are quaternionic and
form the subalgebra $\su(3,\HH)$.  The remaining three rotations ${A_k,
A_{k\ell}, A_{\ell}}$ are the elements of $\gtwo$ that leave invariant the
quaternionic subalgebra spanned by $\{k,\ell,k\ell\}$.  The four
elements~$i,j,i\ell$ and~$j\ell$ can be paired into two complex pairs (of
which~$i+i\ell$, $j+j\ell$ is one choice), and the transformations $\dot A_k$,
$\dot A_{k\ell}$, $\dot A_\ell$ act as~$\su(2,\CC)$ transformations producing
the other pairs of complex numbers.  We henceforth refer to this copy of
$\su(2,\CC)$ as $\suH$.  Hence, the $24$ compact elements form the algebra
$\su(3,\HH)\oplus\suH$, and $\ha=\sla(3,\HH)\oplus\suH$.

\subsection{\boldmath More Subalgebras of $\sla(3,\OO)$}

Consider now the composition $\phiType^* \circ \phiHp^*$.  We define the
subspaces $\ha_1$, $\ha_{23}$, $\ha^\perp_1$, and~$\ha^\perp_{23}$ of
$\sla(3,\OO)$ to be the intersections of pairs of subspaces $\ta_1$,
$\ta_{23}$, $\ha$, and~$\ha^\perp$, as indicated in Table~\ref{IntTH}.  For
example, $\ha_1=\ha\cap\ta_1$ and $\ha^\perp_{23}=\ha^\perp\cap\ta_{23}$.
Table~\ref{IntTH} also indicates the number of basis elements which are boosts
and rotations in each of these spaces, whose commutation rules are given in
Table~\ref{CommTH}.

\begin{table}[tbp]
\begin{center}
 \begin{tabular}{c|cccc|cc}
$\cap$ & $\ta_1$ & $\ta_{23}$ &
	\hspace{2cm}
	& $\cap$ & $\ta_1$ & $\ta_{23}$\\
\cline{1-3}\cline{5-7}
$\ha$ & $\ha_1$ & $\ha_{23}$ & & $\ha$ & $(16,6)$ & $(8,8)$ \\
$\ha^\perp$ & $\ha^\perp_1$ & $\ha^\perp_{23}$ & &
	$\ha^\perp$ & $(20,4)$ & $(8,8)$\\
\multicolumn{3}{c}{Subspace} & & \multicolumn{3}{c}{Signature}
\end{tabular}
\caption{Intersections of Subspaces~$\ta_1$,~$\ta_{23}$,~$\ha$,
and~$\ha^\perp$.}
\label{IntTH}
\end{center}
\end{table}

\begin{table}[tbp]
\begin{center}
\begin{tabular}{|c|c|c|c|c|}
\hline
$\left[ \hspace{.15cm} , \hspace{.15cm} \right]$ &
	$\ha_1$ & $\ha_{23}$ & $\ha^\perp_1$ & $\ha^\perp_{23}$ \\
\hline
$\ha_1$ & $\ha_1$ & $\ha_{23}$ & $\ha^\perp_1$ & $\ha^\perp_{23}$ \\
\hline
$\ha_{23}$ & $\ha_{23}$ & $\ha_1$ & $\ha^\perp_{23}$ & $\ha^\perp_1$ \\ 
\hline
$\ha^\perp_1 $ & $\ha^\perp_1$ & $\ha^\perp_{23}$ & $\ha_1$ & $\ha_{23}$\\
\hline
$\ha^\perp_{23}$ & $\ha^\perp_{23}$ & $\ha^\perp_1$ & $\ha_{23}$ & $\ha_1$ \\
\hline
\end{tabular}
\caption{Commutation structure of~$\ha_1, \ha_{23}, \ha^\perp_1,$
and~$\ha^\perp_{23}$.}
\label{CommTH}.
\end{center}
\end{table}

We see from Table~\ref{CommTH} that, in addition to the subalgebras
$\ha = \ha_1 \oplus \ha_{23}$ and $\ta_1 = \ha_1 \oplus \ha^\perp_1$ constructed
previously, $\ha_1$ and $\ha_1 \oplus \ha^\perp_{23}$ are also subalgebras
of~$\sla(3,\OO)$.  Furthermore, $\phiType^* \circ \phiHp^*$ fixes
\begin{equation}
\phiType^* \circ \phiHp^* \left(\sla(3,\OO)\right) \cap \sla(3,\OO)
  = \ha_1 \oplus \ha^\perp_{23}
\label{H1andH23perp}
\end{equation}
as a subalgebra, since everything in~$\ha_1$ is fixed by both maps, while
everything in~$\ha^\perp_{23}$ is multiplied by~$\xi^2 = -1$.  Thus, even
though the composition of associated Cartan maps is not quite an associated
Cartan map itself (due to the minus sign), it does lead to another
$\ZZ_2$-grading; we return to this point in Section~\ref{cgroup} below.

The subalgebra~$\pa$ in~(\ref{H1andH23perp}) has dimension~$38$ and
signature~$(24,14)$, and is thus a real form of~\hbox{$\aa_5\oplus\aa_1$}.
However, it has a fundamentally different basis
from~$\ha=\sla(3,\HH)\oplus\suH$, as~$\pa$ uses a mixture of the quaternionic
transformations of type~$1$ with the orthogonal-quaternionic transformations
of type~$2$ and type~$3$, while~$\sla(3,\HH)$ is comprised of only the
quaternionic transformations.  We therefore refer to this algebra as
$\sla(2,1,\HH)\oplus\suII$.  Explicitly,
\begin{equation}
\suII
  = \langle \dot G_k + 2\dot S^1_k, \dot G_{k\ell}
	+ 2\dot S^1_{k\ell}, \dot G_\ell + 2\dot S^1_\ell \rangle
\end{equation}
again corresponds to permutations of $\lbrace i,j,j\ell,i\ell \rbrace$ and
fixes $\lbrace k, k\ell, \ell \rbrace$, but is not in~$\gtwo$.  We also note
that the maximal compact subalgebra $\ga'_c$ of
$\phiType^*\circ\phiHp^*\left( \sla(3,\OO)\right)$ has dimension~$36$, and its
pre-image in~$\sla(3,\OO)$ has signature~$(24,12)$.  We identify this
subalgebra as~$\su(3,1,\HH)_2$, since this real algebra has a different basis
than our previously identified~$\su(3,1,\HH)$, which we henceforth refer to
as~$\su(3,1,\HH)_1$.

We summarize in Table~\ref{Maximal} the subalgebras constructed from our three
associated Cartan maps, as well as compositions of these maps.  For each
map~$\phi^*$, the second column lists the signature
of~$\ga'=\phi^*(\sla(3,\OO))$, the third column identifies, and gives the
signature of, the fixed subalgebra $\pa = \phi^*\left(\sla(3,\OO)\right) \cap
\sla(3,\OO)$, and the fourth column identifies the pre-image of the maximal
compact subalgebra of~$\ga' = \phi^*\left(\sla(3,\OO)\right)$, listing both
the algebra and its signature.  Each algebra listed in the third and fourth
columns of Table~\ref{Maximal} is a subalgebra of $\sla(3,\OO)$, constructed
using the other real forms of $e_6$.

\begin{table}[tbp]
{\small
\begin{center}
\begin{tabular}{|cccc|}
\hline
 & Signature of & Signature of & Signature of \\
Map & $\ga' = \phi^*\left(\sla(3,\OO)\right)$ &
	$\pa = \ga' \cap \sla(3,\OO)$ & $(\phi^*)^{-1}(\ga'_c)$ \\
\hline
$1$ & $(52,26)$ & $(52,26)$ & $(52,0)$ \\
& & $\sla(3,\OO)$ & $\su(3,\OO)$ \\
\hline
$\phit^*$ & $(78,0)$ & $(52,0)$ & $(52,26)$ \\
& & $\su(3,\OO)$ & $\sla(3,\OO)$ \\
\hline
$\phiType^*$ & $(52,26)$ & $(36,10)$ & $(36,16)$ \\
& & $\sla(2,\OO) \oplus \uonem$ & $\su(2,1,\OO)$ \\
\hline
$\phiHp^*$& $(36,42)$ & $(24,14)$ & $(24,12)$ \\
& & $\sla(3,\HH)\oplus\suH$ & $\su(3,1,\HH)_1$ \\
\hline
$\phiType^*\circ \phit^*$ & $(46,32)$ & $(36,16)$ & $(36,10)$ \\
& & $\su(2,1,\OO)$ & $\sla(2,\OO) \oplus \uonem$ \\
\hline
$\phiHp^*\circ \phit^*$ & $(38,40)$ & $(24,12)$ & $(24,14)$ \\
& & $\su(3,1,\HH)_1$ & $\sla(3,\HH)\oplus\suH$ \\
\hline
$\phiType^*\circ \phiHp^*$ & $(36,42)$ & $(24,14)$ & $(24,12)$ \\
& & $\sla(2,1,\HH)\oplus\suII$ & $\su(3,1,\HH)_2$ \\
\hline
$\phiType^*\circ \phiHp^*\circ \phit^*$
	& $(38,40)$ & $(24,12)$ & $(24,14)$ \\
& & $\su(3,1,\HH)_2$ & $\sla(2,1,\HH)\oplus\suII$ \\
\hline
\end{tabular}
\caption{Compositions of associated Cartan maps and the corresponding
subalgebras of~$\sla(3,\OO)$.}
\label{Maximal}
\end{center}
}
\end{table}

Despite recognizing~$\ha_1$,~$\ha_1 \oplus \ha^\perp_{23}$,~$\ha_1 \oplus
\ha_{23}$, and~$\ha_1 \oplus \ha^\perp_1$ as subalgebras of~$\sla(3,\OO)$ using
the commutation relations in Table~\ref{CommTH}, we note that~$\su(3,1,\HH)_2$
is not any of these subalgebras.  In the next section, we use another
technique involving associated Cartan maps to give a finer refinement of
subspaces of~$\sla(3,\OO)$, allowing us to provide a nice basis for
$\su(3,1,\HH)_2$ and other subalgebras of~$\sla(3,\OO)$.

\section{\boldmath Constructing Subalgebras of~$\sla(3,\OO)$}
\label{Auto3Apps}

\subsection{Using Composition of Associated Cartan Maps}

We have already used associated Cartan maps to identify the maximal
subalgebras of~$\sla(3,\OO)$.  This was done by separating the algebra into two
separate spaces, one of which was left invariant by the map.  In this section,
we use composition of associated Cartan maps to separate~$\sla(3,\OO)$ into four
or more subspaces spanned by either the compact or non-compact generators,
with the condition that the map either preserves the entire subspace or
changes the character of all the basis elements in the subspace.  We identify
additional subalgebras of~$\sla(3,\OO)$ by taking various combinations of these
subspaces.

We continue to use the associated Cartan maps $\phit^*$, $\phiType^*$,
and~$\phiHp^*$, as well as the subspaces $\ra$, $\ba$, $\ta_1$,
$\ta_{23}$, $\ha$, and~$\ha^\perp$ defined in the previous section.

We first consider the composition $\phiHp^* \circ \phit^*$.  This map fixes
the subspaces \hbox{$\ra_\HH = \ra \cap \ha$}, consisting of quaternionic
rotations, as well as $\ba_\perp = \ba \cap \ha^\perp$, consisting of
orthogonal-quaternionic boosts.  Under $\phiHp^*\circ\phit^*$, the two
subspaces~$\ba_\HH$ and~\hbox{$\ra_\perp = \ra \cap \ha^\perp$} change
signature.  These spaces consist of orthogonal-quaternionic rotations and
quaternionic boosts, respectively.  The dimensions of these four spaces are
displayed in Table~\ref{SubHT}.

\begin{table}[tbp]
\begin{center}
\begin{tabular}{c|cc}
$\cap$ & $\ha$ & $\ha^\perp$ \\
\hline
$\ra$ & $|\ra_\HH| = 24$ & $|\ra_\perp| = 28$ \\
$\ba$ & $|\ba_\HH| = 14$ & $|\ba_\perp| = 12$ \\
\end{tabular}
\caption{Splitting of $\ee_6$ basis under~$\phiHp^* \circ \phit^*$.}
\label{SubHT}
\end{center}
\end{table}

We list the signature of these spaces under $\phiHp^*\circ\phit^*$, and can
thus identify subalgebras of $\phiHp^*\circ\phit^*\left(\sla(3,\OO)\right)$.
However, we are primarily interested in the pre-image of these subalgebras in
our preferred algebra~$\sla(3,\OO)$. Determining the signature of these spaces
in $\sla(3,\OO)$ is straightforward, as the rotations are compact and the
boosts are not.

We use the subspaces represented in Table~\ref{SubHT} to identify subalgebras
of~$\sla(3,\OO)$.  As previously identified, the $24$ rotations in $\ra_\HH$
fixed by the automorphism form the subalgebra~$\su(3,\HH)_1\oplus\suH$, where
$\su(3,\HH)_1\subset\su(3,1,\HH)_1$.  The entries in the first column of
Table~\ref{SubHT} represent all quaternionic rotations and boosts, and form
the subalgebra $\sla(3,\HH)\oplus\suH$.  Of course, the entries in the first
row, $\ra_\HH$ and~$\ra_\perp$, form~$\ff_4$.  We finally consider the entries
$\ra_\HH$ and~$\ba_\perp$ on the main diagonal of the table.  Two
orthogonal-quaternionic boosts commute to a quaternionic rotation, and an
orthogonal-quaternionic boost commuted with a quaternionic rotation is again
an orthogonal-quaternionic boost.  Hence, the subspace
$\ra_\HH\oplus\ba_\perp$ closes under commutation and is a subalgebra with
signature~$(24,12)$.  While both $\so(9-n,n,\RR)$ and~$\su(4-n,n,\HH)$ have
dimension~$36$, only~$\su(3,1,\HH)$ has~$12$ boosts, so $\ra_\HH \oplus
\ba_\perp$ is the previously identified~$\su(3,1,\HH)_1$.

The algebra~$\ra_\HH \oplus \ba_\perp = \su(3,1,\HH)_1$ is a real form of the
complex Lie algebra~$\cc_4$.  Since $\cc_4$ contains $\cc_3$ but not $\bb_3$,
the $21$-dimensional subalgebra contained within $\ra_\HH$ is a real form of
$\cc_3$, not of $\bb_3$.  In addition, any simple $21$-dimensional subalgebra
of $\ra_\HH \oplus \ba_\perp$ is a real from of $\cc_3$.  Eliminating the
boosts from $\su(3,1,\HH)_1$ leaves $\su(3,\HH)_1$.

We next consider the composition $\phiType^* \circ \phit^*$.  This
automorphism separates the basis for~$\sla(3,\OO)$ into the subspaces
\begin{equation}
\begin{array}{ccc}
\ra_1 = \ra \cap \ta_1 &
	\hspace{2cm} &
	 \ba_1 = \ba \cap \ta_1 \\
\ba_{23} = \ba \cap \ta_{23} &
	\hspace{2cm} &
	\ra_{23} = \ra \cap \ta_{23}\\
\end{array}
\end{equation}
As shown in Table~\ref{SubT}, this map leaves the signatures of~$\ra_1$
and~$\ba_{23}$ alone, while it reverses the signatures of~$\ra_{23}$
and~$\ba_1$.

\begin{table}[tbp]
\begin{center}
\begin{tabular}{c|cc}
$\cap$ & $\ta_1$ & $\ta_{23}$ \\
\hline
$\ra$ & $|\ra_1| = 36$ & $|\ra_{23}| = 16$ \\
$\ba$ & $|\ba_1| = 10$ & $|\ba_{23}| = 16$ \\
\end{tabular}
\caption{Splitting of $\ee_6$ basis under $\phiType^* \circ \phit^*$.}
\label{SubT}
\end{center}
\end{table}

We again use the subspaces represented in Table~\ref{SubT} to identify
subalgebras of~$\sla(3,\OO)$.  The subspace~$\ra_1$ is the
subalgebra~$\so(9,\RR)$, containing all subalgebras~$\so(n,\RR)$ for~$n\le9$,
and we have already seen that the subspace $\ta_1=\ra_1\oplus\ba_1$ is
$\so(9,1,\RR)\oplus\uonem$.  Again, the complete set of rotations
$\ra_1\oplus\ra_{23}$ form the subalgebra~$\su(3,\OO)$, which is a real form
of~$\ff_4$.  Interestingly, the subspace~$\ra_1 \oplus \ba_{23}$ on the main
diagonal is another form of~$\ff_4$.  The $16$ boosts in~$\ra_1\oplus\ba_{23}$
identify this form of~$\ff_4$ as~$\su(2,1,\OO)$.

We finally consider the composition $\phiType^*\circ\phiHp^*\circ\phit^*$,
which creates a finer refinement than compositions of two maps.  The resulting
subspaces are listed in Table~\ref{SubH}.  We continue with our previous
conventions for designating intersections of subspaces, that is,
$\ra_{23,\perp} = \ra \cap \ta_{23} \cap \ha^\perp$.

\begin{table}[tbp]
\begin{center}
\begin{tabular}{c|cccc}
$\cap$ & $\ha_1$ & $\ha_{23}$ & $\ha^\perp_1$ & $\ha^\perp_{23}$ \\
\hline
$\ra$ & $|\ra_{1,\HH}|=16$ & $|\ra_{23,\HH}|=8$ &
	$|\ra_{1,\perp}|=20$ & $|\ra_{23,\perp}|=8$ \\
$\ba$ & $|\ba_{1,\HH}|=6$ & $|\ba_{23,\HH}|=8$ &
	$|\ba_{1,\perp}|=4$ & $|\ba_{23,\perp}|=8$ \\
\end{tabular}
\caption{Splitting of~$\ee_6$ basis under
$\phiType^*\circ\phiHp^*\circ\phit^*$.}
\label{SubH}
\end{center}
\end{table}

Using this division of~$\sla(3,\OO)$, we find a large list of subalgebras
of~$\sla(3,\OO)$ simply by combining certain subspaces.  The subspace
description of these algebras, as well as their identity and signature
in~$\sla(3,\OO)$, is listed in Table~\ref{Refine}.  This fine refinement
of~$\sla(3,\OO)$ provides a description of the basis for~$\su(3,1,\HH)_1$
and~$\su(3,1,\HH)_2$, as well as~$\sla(3,\HH)$ and~$\sla(2,1,\HH)$.
%

\begin{table}[tbp]
\begin{center}
\begin{tabular}{|l|l|c|c|}
\hline
Basis & Subalgebra of $\sla(3,\OO)$ & Signature \\
\hline
$\ra_{1,\HH}$ & $\su(2,\HH)\oplus\suH \oplus \su(2)$ & $(10+3+3,0)$ \\
\hline
$\ra_{1,\HH}\oplus \ra_{23,\perp}$ & $\su(3,\HH)_2\oplus\suH$ & $(21+3,0)$ \\
\hline
$\ha_1=\ra_{1,\HH}\oplus \ba_{1,\HH}$ & $\sla(2,\HH)\oplus\suH $ &
	$(10+3+3,5+1)$ \\
  & \hspace{.5cm} $\oplus \su(2) \oplus \uonem$ & \\
\hline
$\ra_{1,\HH}\oplus \ba_{23,\perp}$ & $\su(2,1,\HH)_1\oplus\suH$ & $(13+3,8)$ \\
\hline
$\ra_{1,\HH}\oplus \ba_{23,\HH}$ & $\su(2,1,\HH)_2\oplus\suH$ & $(13+3,8)$ \\
\hline
$\ra_\HH=\ra_{1,\HH}\oplus \ra_{23,\HH}$ & $\su(3,\HH)\oplus\suH$ & $(21+3,0)$ \\
\hline
$\ra_{1,\HH}\oplus \ba_{1,\perp}$ & $\so(5,\RR)\oplus\so(4,1,\RR)$ & $(10+6,4)$ \\
\hline
$\ra_1=\ra_{1,\HH}\oplus \ra_{1,\perp}$ & $\so(9) = \su(2,\OO)$ & $(36,0)$ \\
\hline
$\ha=\ra_{1,\HH}\oplus \ba_{1,\HH}$ & $\sla(3,\HH) \oplus\suH $ & $(21+3,14)$ \\
  \hspace{.5cm} $\oplus \ba_{23,\HH}\oplus \ra_{23,\HH}$ & & \\
\hline
$\ha_1\oplus\ha^\perp_{23}=\ra_{1,\HH}\oplus \ba_{1,\HH}$ &
	$\sla(2,1,\HH)_1\oplus\suII$ & $(21+3,14)$ \\
  \hspace{2cm} $\oplus \ra_{23,\perp}\oplus \ba_{23,\perp}$ & & \\
\hline
$\ta_1=\ra_{1,\HH}\oplus \ba_{1,\HH}$ & $\sla(2,\OO)\oplus \uonem$ & $(36,9+1)$ \\
  \hspace{.5cm} $ \oplus \ba_{1,\perp}\oplus \ra_{1,\perp}$ & & \\
\hline
$\ra_\HH\oplus\ba_\perp=\ra_{1,\HH}\oplus \ba_{23,\perp}$ &
	$\su(3,1,\HH)_1$ & $(24,12)$ \\
  \hspace{2cm}$\oplus \ra_{23,\HH}\oplus \ba_{1,\perp}$ & & \\
\hline
$\ra_{1,\HH}\oplus \ra_{23,\perp}$ & $\su(3,1,\HH)_2$ & $(24,12)$ \\
  \hspace{.5cm} $\oplus \ba_{23,\HH}\oplus \ba_{1,\perp}$ & & \\
\hline
$\ra=\ra_{1,\HH}\oplus \ra_{23,\perp}$ & $\su(3,\OO)$ & $(52,0)$ \\
  \hspace{.5cm} $\oplus \ra_{23,\HH}\oplus \ra_{1,\perp}$ & & \\
\hline
$\ra_1\oplus\ba_{23}=\ra_{1,\HH}\oplus \ba_{23,\perp}$ &
	$\su(2,1,\OO)$ & $(36,16)$ \\
  \hspace{2cm} $\oplus \ba_{23,\HH}\oplus \ra_{1,\perp}$ & & \\
\hline
\end{tabular}
\caption{Subalgebras of~$\sla(3,\OO)$ using Cartan decompositions.}
\label{Refine}
\end{center}
\end{table}

\subsection{\boldmath Chains of Subalgebras of~$\sla(3,\OO)$}
\label{Chains}

We have used associated Cartan maps to produce large simple subalgebras
of~$\sla(3,\OO)$, ranging in dimension from~$52$, for~$\ff_4$, to~$21$,
for~$\cc_3$.  Each of these subalgebras in turn has its own associated Cartan
maps, which we could use to find even smaller subalgebras, thereby giving a
catalog of subalgebra chains contained within~$\sla(3,\OO)$.  However, having
identified the real form of the large subalgebras of~$\sla(3,\OO)$, it is not
too difficult a task to find the smaller algebras simply by looking for simple
subalgebras of smaller dimension and/or rank, using the tables of real forms
listed in~\cite{gilmore} when needed.  Furthermore, we can choose our smaller
subalgebras and their bases so that they use a subset of our preferred basis
for the Cartan subalgebra of $\sla(3,\OO)$, namely
$\lbrace \dot B^1_{tz}, \dotBtz, \dot R^1_{x\ell},
 \dot S^1_{\ell}, \dot G_{\ell}, \dot A_{\ell} \rbrace$,
henceforth referred to as our \textit{Cartan basis}.
%

\goodbreak

\begin{figure}[tbp]
\null\vspace{-0.25in}
\hfill
\xymatrixcolsep{5pt}
\hspace{-1in}
\xymatrix@M=3pt@H=10pt{
 & & & *++\txt{$\sla(3,\OO)$ \\ $[\ee_6]$} & \\
\save[]+<0cm,.4cm>
*++\txt{$\sla(2,1,\HH) \oplus\suII$ \\ $[\aa_5 \oplus \aa_1]$}
  \ar@{^{(}.>}[urrr]<2.75ex>_{} \restore & & & &
\save[]+<0cm,.4cm>
*++\txt{$\sla(3,\HH) \oplus\suH$ \\ $[\aa_5 \oplus \aa_1]$}
\ar@{_{(}.>}[ul]_{} \restore \\
*++\txt{$\sla(2,1,\HH)$ \\ $[\aa_5]$} \ar@{.>}[u]^{G_\ell+2S^1_\ell} &
*++\txt{$\su(2,1,\OO)$ \\ $[\ff_{4(36,16)}]$}
  \ar@/^3pc/[uurr]^(.3){B^1_{tz}}^(.24){\Btz} &
*++\txt{$\su(3,\OO)$ \\ $[\ff_4]$} \ar[uur]^(.4){B^1_{tz}}^(.3){\Btz} &
*++\txt{$\sla(2,\OO) = \so(9,1)$ \\ $[\dd_5]$} \ar@{.>}[uu]_{\Btz} &
*++\txt{$\sla(3,\HH)$ \\ $[\aa_5]$}  \ar@{.>}[u]_{A_\ell} \\ & & & & \\
 & *++\txt{$\su(3,1,\HH)_1$ \\ $[\cc_4]$}
  \ar@/^4pc/[uuuurr]^(.7){B^1_{tz}}^(.64){\Btz}
 & *++\txt{$\su(2,\OO) = \so(9)$ \\ $[\bb_4]$}
  \ar@{_{(}.>}[uul]^(.4){} \ar@{_{(}.>}[uu]^(.4){} \ar[uur]^(.65){B^1_{tz}} &
*++\txt{$\so(8)$ \\ $[\dd_4]$} \ar@{_{(}.>}[l]^(.4){} & \\ & & & & \\
*++\txt{$\su(2,1,\HH)_1$ \\ $[\cc_3]$}
  \ar[uuuu]<1ex>^(.55){B^1_{tz}}^(.45){\Btz}
  \ar@{.>}[uuuur]^(.5){G_\ell+2S^1_\ell}
  \ar@{.>}[uur]{} \ar[uur]<1ex>_(.35){~G_\ell+2S^1_\ell} &
*++\txt{$\su(3,\HH)_1$ \\ $[\cc_3]$}
  \ar `u[r] `[rrru] [uurrruu] ^(.6){B^1_{tz}} ^(.5){\Btz}
  \ar@{.>}[uu]<2ex>{} \ar[uu]<1ex>_(.5){A_\ell} \ar@{.>}[uuuur]^(.5){A_\ell} &
*++\txt{$\su(1,\OO) = \so(7)$ \\ $[\bb_3]$}
  \ar[uur]^(.5){R^1_{x\ell}} \ar@{.>}[uu]_(.5){R^1_{x\ell}} &
*++\txt{$\sla(2,\HH)$ \\ $[\aa_3 = \dd_3]$}
  \ar `l[uuu] `[uuull] `[uuulll]_(.9){A_\ell, \Btz} [uuuulll]
  \ar@/_4pc/[uuuu]^(.7){G_\ell+2S^1_\ell}^(.8){A_\ell}
  \ar@/_8pc/[uuuur]^(.7){G_\ell+2S^1_\ell}^(.6){\Btz} & \\ & & & & \\
 & *++\txt{$\su(2,\HH) = sp(2) $ \\ $[\bb_2 = \cc_2]$}
  \ar[uuuur]^(.8){A_\ell}^(.72){G_1+2S^1_\ell}
  \ar@/_1.8pc/[uurr]^(.3){B^1_{tz}}
  \ar@{.>}[uul]^(.5){A_\ell} \ar@{.>}[uu]<1ex>^(.65){G_\ell+2S^1_\ell} &
*++\txt{$[\gtwo]$} &
*++\txt{$\su(3,\CC)_s$ \\ $[\aa_2]$} \ar@{.>}[uull]_(.6){G_\ell} &
*++\txt{$\sla(2,\CC)_s$ \\ $[\aa_1\oplus \aa_1]$ }
  \ar@{.>}[uul]<1ex>{} \ar[uul]_(.5){G_\ell-S^1_\ell} \\
 & & *++\txt{$\so(4,\RR)$ \\ $[\dd_2 = \cc_1 \oplus \cc_1]$}
  \ar@{_{(}.>}[ul]^(.4){} \ar[uuuuur]^(.6){A_\ell}^(.55){G_\ell+2S^1_\ell} & & \\
 & & & *++\txt{$\su(1,\HH)$ \\ $[\cc_1]$}
  \ar[uuuul]<-1ex>_(.28){A_\ell}_(.2){G_\ell+2S^1_\ell}
  \ar@{.>}[ul]<1ex>{} \ar[ul]_(.45){R^1_{x\ell}} &
*++\txt{$\su(2,\CC)_s$ \\ $[\aa_1]$}
  \ar[uu]_(.5){B^1_{tz}} \ar@{.>}[uul]<1ex>^(.5){}
  \ar[uul]_(.5){R^2_{x\ell} \to S^1_\ell } \\
 & & & \uone \ar[u]_(.4){G_\ell - S^1_\ell} & \uone \ar[u]_(.4){R^1_{x\ell}}\\
}
\caption{Preferred subalgebra chains of~$\ee_6$.}
\label{SubChain}
\end{figure}
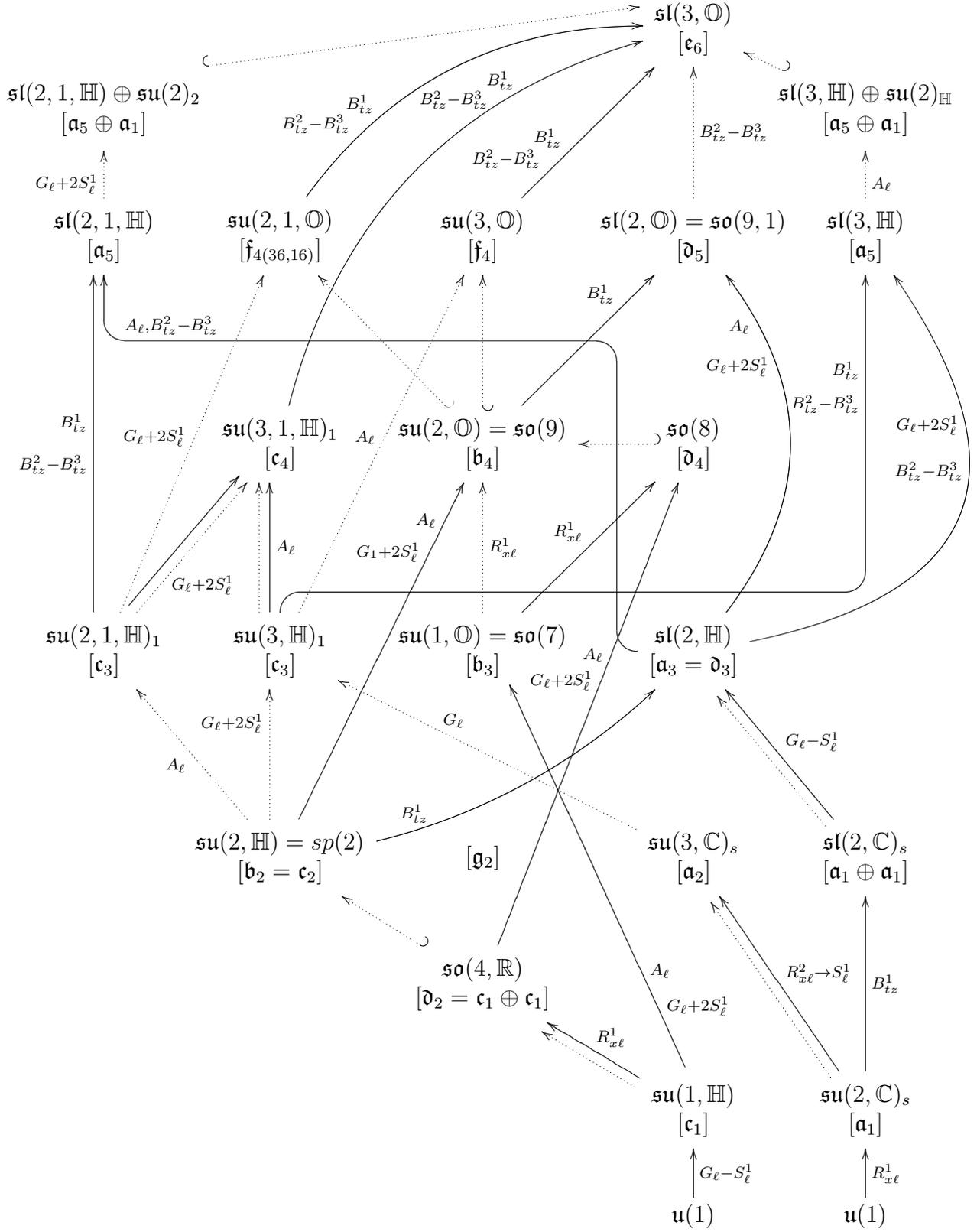

We display these chains of subalgebras of~$\sla(3,\OO)$ in the following
tables.  Each table is built from a $\uone$ algebra, generated by a single
element of our Cartan basis.  We extend each algebra~$\ga$ to a larger
algebra~$\ga'$ by adding elements to the basis for~$\ga$.  In particular, each
algebra of higher rank must add new elements of our Cartan basis, as indicated
along the arrows (with dots suppressed).  Figure~\ref{SubChain} is built from
the algebra \hbox{$\uone = \langle G_l - S^1_l \rangle$}, leading to
\begin{equation}
\uone \subset \su(1,\HH) \subset \su(2,\HH) \subset \su(3,\HH)_1
	\subset \sla(3,\HH)
\end{equation}
Additional subalgebras of~$\sla(3,\OO)$ can be inserted into this chain
of subalgebras.  For instance, we can insert~$\sla(2,\HH)$ between~$\su(2,\HH)$
and $\sla(3,\HH)$, and extend~$\sla(2,\HH)$ to~$\sla(2,\OO)$.  We can also expand
$\su(1,\HH)$ to \hbox{$\su(1,\OO) = \so(7)$} and insert the~$\so(n,\RR)$ chain
for~$n \ge 7$ into the figure.  However, $\so(7)$ uses a basis in this chain
that is not compatible with \hbox{$\gtwo = \hbox{aut}(\OO)$}.  We do not list
all of the possible $\uone$ subalgebras, but do add the subalgebra 
\hbox{$\uone = \langle R^1_{x\ell} \rangle$} into the chain and extend it to
$\su(2,\CC)_s$ and $\sla(2,\CC)_s$, which use the standard (type~1) matrix
definition of $\su(2,\CC)$ and~$\sla(2,\CC)$.

With one exception, the algebras in Figure~\ref{SubChain} are built from
subalgebras by extending the basis at each step; we do not allow any changes
of basis.  The one exception is the inclusion of the chain
\begin{equation}
\su(2,\CC)_s \subset \su(3,\CC)_s \subset \su(3,\HH)_1
\end{equation}
using the Cartan basis elements $\dot{S}^1_\ell$ and $\dot{G}_\ell$.
We also use different notations to indicate possible methods used to identify
the subalgebras, with dashed and solid arrows indicating that the root diagram
of the smaller algebra can be obtained as a \textit{slice} or as a
\textit{projection} of that of the larger algebra; for details,
see~\cite{aaron_thesis}.

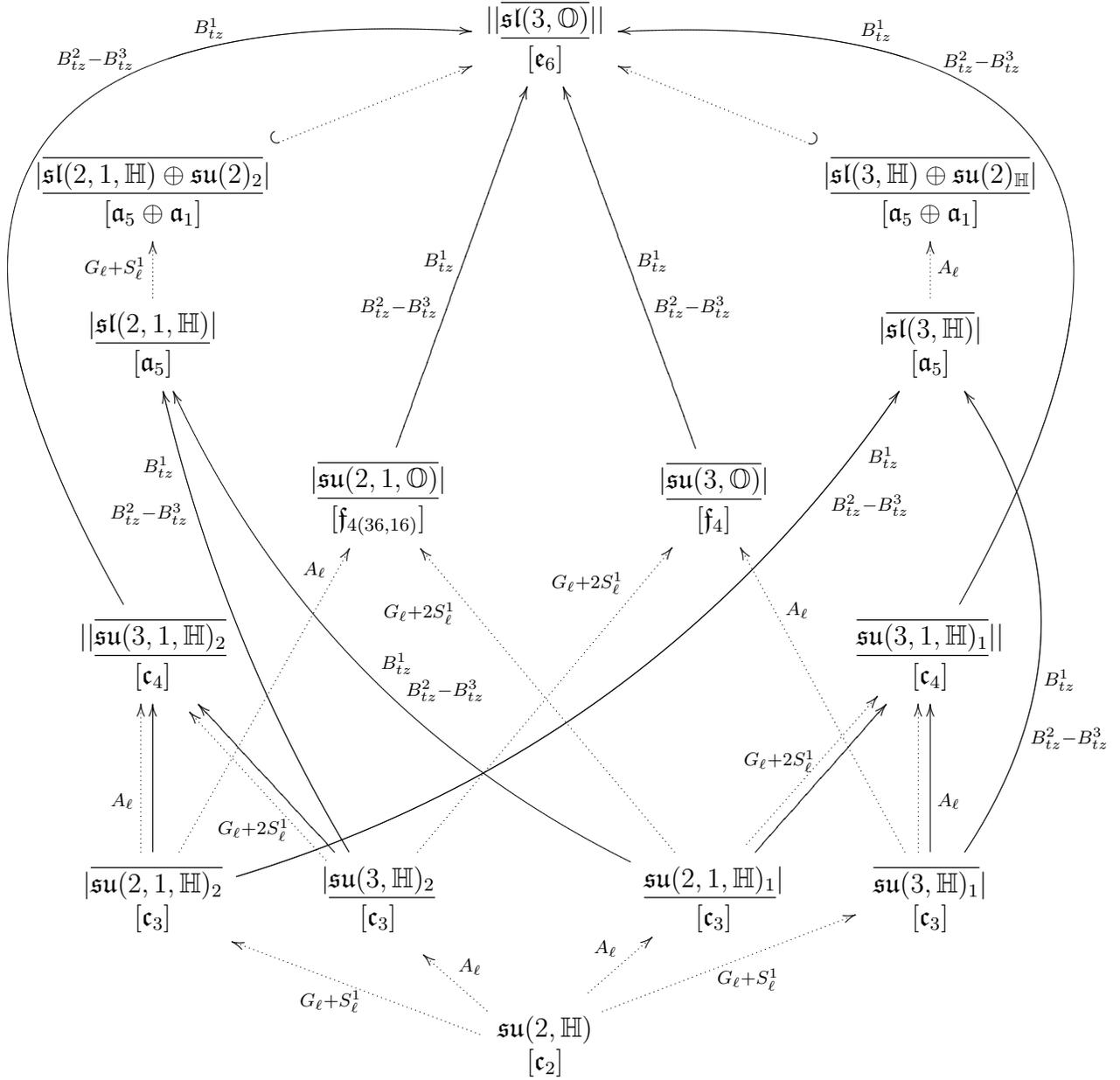
\begin{figure}[tbp]
\centering
\xymatrixcolsep{5pt}
\xymatrix@M=3pt@H=10pt{
 & & *++\txt{$\bbline{\ppline{\rgline{\sla(3,\OO)}}}$ \\ $[\ee_6]$} & & \\
*++\txt{$\bpline{\rgline{\sla(2,1,\HH)\oplus\suII}}$ \\ $[\aa_5 \oplus \aa_1]$}
  \ar@{^{(}.>}[urr]^(.5){} & & & &
*++\txt{$\bpline{\rgline{\sla(3,\HH) \oplus\suH}}$ \\ $[\aa_5 \oplus \aa_1]$}
  \ar@{_{(}.>}[ull]^(.5){} \\
*++\txt{$\bpline{\gline{\sla(2,1,\HH)}}$ \\ $[\aa_5]$}
  \ar@{.>}[u]^(.5){G_\ell + S^1_\ell} & & & &
*++\txt{$\bpline{\rline{\sla(3,\HH)}}$ \\ $[\aa_5]$} \ar@{.>}[u]_(.5){A_\ell} \\
 & *++\txt{$\bpline{\rgline{\su(2,1,\OO)}}$ \\ $[\ff_{4(36,16)}]$}
  \ar[uuur]^(.5){B^1_{tz}}^(.4){\Btz} & &
*++\txt{$\bpline{\rgline{\su(3,\OO)}}$ \\ $[\ff_4]$}
  \ar[uuul]_(.5){B^1_{tz}}_(.4){\Btz} & \\
*++\txt{$\ppline{\rgline{\su(3,1,\HH)_2}}$ \\ $[\cc_4]$}
  \ar@/^13.2pc/[uuuurr]^(.75){B^1_{tz}}^(.65){\Btz~} & & & &
*++\txt{$\bbline{\rgline{\su(3,1,\HH)_1}}$ \\ $[\cc_4]$}
  \ar@/_13pc/[uuuull]_(.75){B^1_{tz}}_(.65){\Btz} \\ & & & & & \\
*++\txt{$\pline{\rline{\su(2,1,\HH)_2}}$ \\ $[\cc_3]$}
  \ar@{.>}[uu]<1ex>^(.4){A_\ell} \ar[uu]^(.4){}
  \ar@{.>}[uuur]^(.8){A_\ell}
  \ar@/_3.5pc/[uuuurrrr]_(.88){B^1_{tz}}_(.82){\Btz} &
*++\txt{$\pline{\gline{\su(3,\HH)_2}}$ \\ $[\cc_3]$}
  \ar@{.>}[uul]<1ex>^(.4){} \ar[uul]^(.35){G_\ell+2S^1_\ell\,}
  \ar@{.>}[uuurr]^(.75){G_\ell+2S^1_\ell}
  \ar@/^1pc/[uuuul]^(.80){B^1_{tz}}^(.72){\Btz} & &
*++\txt{$\bline{\gline{\su(2,1,\HH)_1}}$ \\ $[\cc_3]$}
  \ar@{.>}[uur]<1ex>_(.5){} \ar[uur]^(.5){G_\ell+2S^1_\ell}
  \ar@{.>}[uuull]^(.75){G_\ell+2S^1_\ell}
  \ar@/^3pc/[uuuulll]_(.5){B^1_{tz}}_(.45){\Btz}
& *++\txt{$\bline{\rline{\su(3,\HH)_1}}$ \\ $[\cc_3]$}
  \ar@{.>}[uu]<1ex>_(.5){} \ar[uu]_(.4){A_\ell} \ar@{.>}[uuul]_(.7){A_\ell}
  \ar@/_4pc/[uuuu]_(.4){B^1_{tz}}_(.3){\Btz} \\ & &     
*++\txt{$\su(2,\HH)$ \\ $[\cc_2]$}
  \ar@{.>}[ull]<1ex>^(.5){G_\ell + S^1_\ell} \ar@{.>}[ul]<1ex>_(.5){A_\ell}
  \ar@{.>}[ur]<1ex>^(.5){A_\ell} \ar@{.>}[urr]<1ex>_(.5){G_\ell + S^1_\ell}
}
\caption{Four real forms of~$\cc_3$.  An underline indicates the presence of
$\suII$, an overline indicates the presence of $\suH$, each initial vertical
line indicates the presence of one component of $\ba_{23,\HH}\oplus
\ra_{23,\perp}$, and each final vertical line indicates the presence of one
component of $\ra_{23,\HH} \oplus \ba_{23,\perp}$.}
\label{C3algs}
\end{figure}

Finally, we have identified four different real forms of~$\cc_3$, all of which
contain~$\su(2,\HH)$.  Space constraints limit us to listing only
$\su(2,1,\HH)_1$ and $\su(3,\HH)_1$ in Figure~\ref{SubChain}, but the algebras
$\su(2,1,\HH)_2$, $\su(3,\HH)_2$, and $\su(3,1,\HH)_2$ should also be in this
table.  We list these four real forms of~$\cc_3$ algebras, all built
from~$\su(2,\HH)$, in Figure~\ref{C3algs}, and include all the algebras which
are built from the~$\cc_3$ algebras.  Figure~\ref{C3algs} can be incorporated
into Figure~\ref{SubChain} without having to adjust our choice of Cartan
basis.

\begin{figure}[tbp]
{\small
\begin{center}
\begin{minipage}{6in}
\begin{center}
\xymatrixcolsep{5pt}
\xymatrix@M=6pt@H=48pt@R=48pt{
 & *++\txt{$\ee_{6 (78,0)} $ \\
  $\arrowvert \phit(\sla(3,\OO)) \arrowvert = 78$} & \\
*++\txt{$\ee_{6 (52,26)}$ \\
	$\arrowvert\phit(\su(3,\OO))
	\arrowvert = 52$} \ar@/^3pc/@{<->}[ur]^(.35){\phit} &
*++\txt{$\ee_{6 (46,32)}$ \\
	$\arrowvert\phit(\so(9,1) \oplus \uonem)\arrowvert = 46$}
	\ar@{<->}[u]^(.5){\phiType} &
*++\txt{$\ee_{6 (38,40)}$ \\
	$\arrowvert\phit(\sla(3,\HH)\oplus\suH)
	\arrowvert = 38$} \ar@/_3pc/@{<->}[ul]_(.35){\phiHp} \\
*++\txt{$\ee_{6 (52,26)}$ \\
	$\arrowvert\phit(\su(2,1,\OO))\arrowvert = 52$}
	\ar@/^1pc/@{<->}[ur]^(.35){\phit} \ar@{<->}[u]^(.5){\phiType} &
*++\txt{$\ee_{6 (36,42)}$ \\
	$\arrowvert\phit(\su(2,1,\HH))\arrowvert = 36$}
	\ar@/_1pc/@{<->}[ul]_(.35){\phiHp}
	\ar@/^1pc/@{<->}[ur]^(.35){\phit} &
*++\txt{$\ee_{6 (38,40)}$ \\
	$\arrowvert\phit(\sla(2,1, \HH)\oplus\suII)
	\arrowvert = 38$}
	\ar@/_1pc/@{<->}[ul]_(.35){\phiHp}
	\ar@{<->}[u]_(.5){\phiType} \\
 & *++\txt{$\ee_{6 (36,42)}$ \\
	$\arrowvert\phit(\su(3,1,\HH)_2)\arrowvert = 36$}
	\ar@/^3pc/@{<->}[ul]^(.65){\phiHp}
	\ar@{<->}[u]^(.5){\phiType} \ar@/_3pc/@{<->}[ur]_(.65){\phit} & \\
}
\caption{Composition of associated Cartan maps of $\ee_6$ acting on real forms
of $\ee_6$, showing the maximal compact subalgebra under $\phit$.}
\label{AutoII}
\end{center}
\end{minipage}
\end{center}
}
\end{figure}
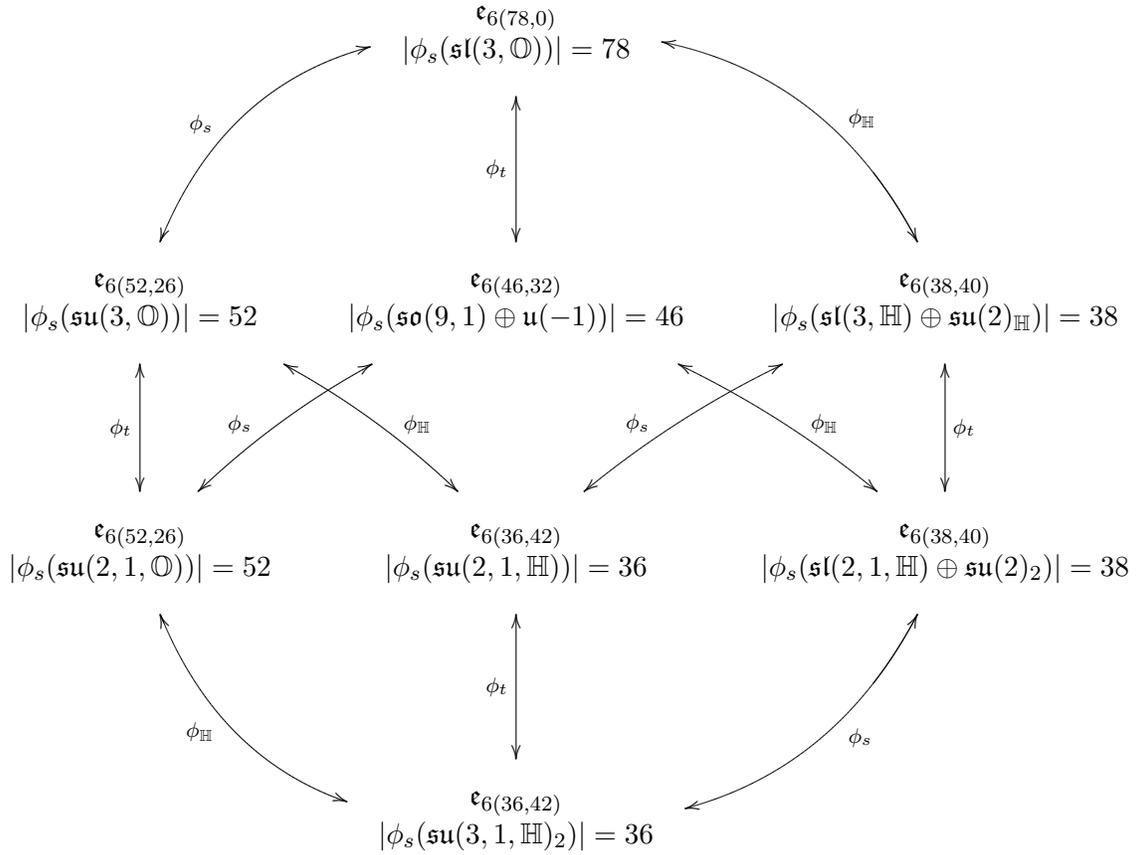

\section{Discussion}
\label{discussion}

\subsection{The Group of Associated Cartan Maps}
\label{cgroup}

We return to the structure of associated Cartan maps~(\ref{associated}).  The
square of~$\phi^*$ is clearly the original involution $\phi$, and the inverse
of $\phi^*$ is obtained by replacing $\xi$ with $-\xi$, or equivalently as the
cube of $\phi^*$.  But the composition of two associated Cartan maps is not
quite an associated Cartan map as we have defined them above, although it does
lead to a graded Lie algebra structure, which we use to define a group
operation as follows.

Given two associated Cartan maps, written symbolically as
\begin{eqnarray}
\phi^*_1(\pa_1+\ma_1) &=&  \pa_1 + \xi \ma_1 \\
\phi^*_2(\pa_2+\ma_2) &=&  \pa_2 + \xi \ma_2
\end{eqnarray}
where
\begin{equation}
\pa_1\oplus\ma_1 = \ga^\CC = \pa_2\oplus\ma_2
\end{equation}
we define their product to be
\begin{equation}
(\phi^*_1 \star \phi^*_2)(q_{pp}+q_{pm}+q_{mp}+q_{mm})
  = q_{pp}+\xi q_{pm}+\xi q_{mp}+q_{mm}
\end{equation}
where
\begin{equation}
q_{pp}\in\pa_1\cap\pa_2
  \qquad
q_{pm}\in\pa_1\cap\ma_2
  \qquad
q_{mp}\in\ma_1\cap\pa_2
  \qquad
q_{mm}\in\ma_1\cap\ma_2
\end{equation}
and which differs from composition by the sign of the last term.  It is easily
verified that
\begin{equation}
\ga^\CC = (q_{pp}\oplus q_{mm}) \oplus (q_{pm}\oplus q_{mp})
\end{equation}
is a $\ZZ_2$-grading of $\ga^\CC$, so that $\phi^*_1\star\phi^*_2$ is indeed
an associated Cartan map.  The operation $\star$ is commutative, and the
square of any associated Cartan map using $\star$ is the identity map.

The set of associated Cartan maps therefore forms a group under the operation
$\star$.  We consider in particular the group generated by the associated
Cartan maps $\phit^*$, $\phiType^*$, and $\phiHp^*$, which is easily seen to
be a copy of $(\ZZ_2)^3$, and hence of order 8.  The orbit of
\hbox{$\sla(3,\OO)=\ee_{6(52,26)}$} under this group is shown in
Figure~\ref{AutoII}, from which the multiplication table can be inferred.

\begin{figure}[tbp]
\centering
\includegraphics[height=2.9in,bb=56 217 548 679]{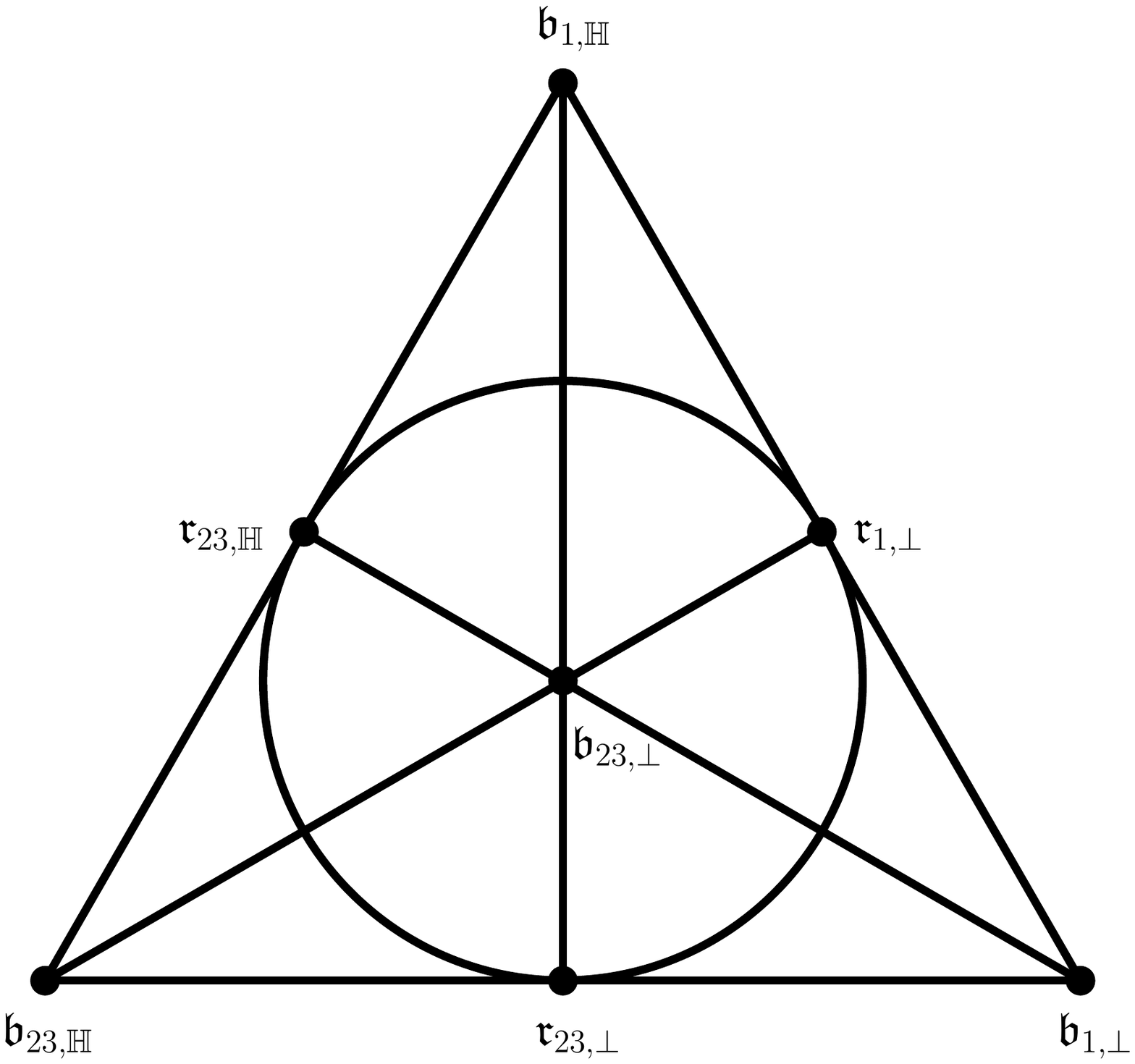}
\caption{The $(\ZZ_2)^3$ multiplication table under commutation for the
subspaces of $\ee_6$ determined by the associated Cartan maps $\phiType^*$,
$\phiHp^*$, and $\phit^*$.  All elements square to $\ra_{1,\HH}$.}
\label{e6fano}
\end{figure}

The multiplication table of the finite group $(\ZZ_2)^3$ is identical to that
of the octonionic units, without the minus signs, and can therefore also be
represented using the 7-point projective plane.  Using commutation as the
operation, the same multiplication table applies directly to the 8 subspaces
listed in Table~\ref{SubH}, as shown in Figure~\ref{e6fano}.  Using octonionic
language for this multiplication table, the 15 proper subalgebras of
$sl(3,\OO)$ listed in Table~\ref{Refine} consist precisely of 1 ``real''
subalgebra, corresponding to the ``identity element'' $\ra_{1,\HH}$, 7
``complex'' subalgebras, formed by adding any one other subspace,
corresponding to the ``points'' in the multiplication table, and 7
``quaternionic'' subalgebras, formed by adding any one additional subspace
(and ensuring the algebra closes), corresponding to the ``lines'' in the
multiplication table.  There is of course one subalgebra missing from this
description, namely the ``octonionic'' algebra $sl(3,\OO)$ itself,
corresponding to the entire Fano ``plane''.

\subsection{\boldmath Real Forms of $\ee_6$}
\label{RealForms}

There are 5 real forms of $\ee_6$, all of which appear in Figure~\ref{AutoII},
although it is at first sight somewhat surprising that several of them appear
more than once.  However, our interpretation of $SL(3,\OO)$ is tied to a
particular choice of basis, so different copies of a given real form yield
different decompositions of $\sla(3,\OO)$, not necessarily with the same
signature.

We regard Figure~\ref{AutoII} itself as an indication that there are really 8
``real forms'' of $\ee_6$ of relevance to the structure of $SL(3,\OO)$, and
hence of possible relevance to physics.  At the very least, all 8 of these
real forms played a role in the construction of the ``maps'' of $\sla(3,\OO)$
given in Figures~\ref{SubChain} and~\ref{C3algs}


\newpage

\section*{Acknowledgments}

This paper is a revised version of Chapter~5 of a dissertation submitted by AW
in partial fulfillment of the degree requirements for his Ph.D.\ in
Mathematics at Oregon State University~\cite{aaron_thesis}.  The revision was
made possible in part through the support of a grant from the John Templeton
Foundation.

\bibliographystyle{unsrt}
\bibliography{e6sub}

\end{document}